\documentclass[runningheads]{llncs}
\usepackage[T1]{fontenc}
\usepackage{graphicx}
\usepackage{amsfonts}
\usepackage{amssymb}
\usepackage{amsbsy}
\usepackage{amsmath}
\usepackage{color}
\usepackage[capitalise, noabbrev]{cleveref}
\usepackage{subcaption} 
\usepackage{pgfplots}

\usepackage[
algo2e,
lined,
boxed,
commentsnumbered,
linesnumbered,
ruled
]{algorithm2e}

\usepackage{pgfkeys}
\def\addlegendimage{\csname pgfplots@addlegendimage\endcsname}

\newcommand{\R}{\mathbb{R}}
\newcommand{\N}{\mathbb{N}}

\newcommand{\J}{\mathcal{J}}
\newcommand{\Jhat}{\hat{\mathcal{J}}}

\newcommand{\cJhatn}{{{\Jhat_\red}}}

\newcommand{\pr}{\textnormal{pr}}
\newcommand{\du}{\textnormal{du}}

\newcommand{\Params}{\mathcal{P}}
\newcommand{\Proj}{\mathrm{P}}

\newcommand{\red}{{N}}

\newcommand{\bdg}{a_{\rm DG}}

\newcommand{\rhs}{l_{\rm DG}}
\newcommand{\n}[0]{n}
\newcommand{\faces}{\mathcal{F}}
\newcommand{\face}{e}

\newcommand{\Grid}{\mathcal{T}_H}
\newcommand{\grid}{\tau_h}
\newcommand{\redspace}{V}
\newcommand{\fullredspace}{V}
\newcommand{\vspan}{\text{ span }}

\newcommand{\FOM}{DG-MsFEM}

\newcommand\restrict[1]{\raisebox{-.5ex}{$\big|$}_{#1}}

\def\V{V}

\begin{document}
\title{Adaptive Localized Reduced Basis Methods for Large Scale PDE-constrained Optimization}
\titlerunning{Localized RB Methods for Large Scale PDE-constrained Optimization}

\author{Tim Keil \and
Mario Ohlberger %\orcidID{0000-0002-6260-3574} 
\and
Felix Schindler
}
\authorrunning{T. Keil \and M. Ohlberger \and F. Schindler}
\institute{Institute for Analysis and Numerics, University of Münster, Einsteinstr. 62, 48149 M\"unster, Germany.
\url{https://www.wwu.de/AMM/ohlberger} 
\email{\{tim.keil,mario.ohlberger,felix.schindler\}@uni-muenster.de}}

\maketitle              % typeset the header of the contribution
\begin{abstract}
In this contribution, we introduce and numerically evaluate a certified and adaptive localized reduced basis method as a local model in a trust-region optimization method for parameter optimization constrained by partial differential equations. 

\end{abstract}
\section{Introduction} \vspace*{-.5em}
We are concerned with efficient and certified approximations of multiscale or large-scale PDE-constrained parameter optimization problems.
In particular, for a parameter space $\Params \subset \R^P$, $P \in \N$, we search for a local solution of \eqref{P}, i.e.\\[-.5em]
{\color{white}
	\begin{equation}
	\tag{P}
	\label{P}
	\end{equation}
}\vspace{-40 pt} % <= adjust this manually to match the document class and font size!
\begin{subequations}\begin{align}
	& \min_{\mu \in \Params} \J(u_{\mu}, \mu),
	&&\text{with } \J: V \times \Params \to \mathbb{R},
	\tag{P.a}\label{P.argmin}\\[-1em]
	\intertext{%
		where $u_{\mu} \in V$ solves the parameterized state equation for $\mu \in \Params$:
	}\nonumber\\[-2em]
	&a(u_\mu, v;\mu) = l(v;\mu) &&\text{for all } v \in V.
	\tag{P.b}\label{P.state}
	\end{align}\end{subequations}%
% ... and reset this counter.
\setcounter{equation}{0}% we are done cheating here
Here, $V$ is a Hilbert space and, for each admissible parameter $\mu \in \Params$, $a(\cdot,\cdot;\mu): V \times V \to \R$ denotes a continuous and coercive bilinear form and $l(\cdot;\mu) \in V'$.

Particularly, we are interested in multiscale or large-scale applications in the sense that the 
state equation (\ref{P.state}) is a weak formulation of a PDE of the form 
\begin{equation}
  \label{eq:o2}
        - \nabla \cdot \big( A(\mu) \nabla u_{\mu} \big) = f(\mu) \ \ \text{in}~ \Omega, \qquad
      u_{\mu} = 0  \ \ \text{on}~ \partial \Omega,
 \end{equation}
where $\Omega \subset \R^d$ is a bounded domain and $A$ denotes a diffusion tensor with a rich structure that would lead to very high dimensional approximation spaces for the state space when approximated with classical finite element type methods.  
On the other hand, employing model order reduction for the involved equations is problematic for high dimensional parameter spaces $\Params$ due to large offline times.

Efficient adaptive model order reduction in combination with trust-region optimization methods in this context has recently been addressed in \cite{MR4269464,MR3716566,WenZhar2022} employing global reduced basis models and in \cite{KeilOhlberger22} employing a reduced basis localized orthogonal decomposition.
The latter work is a first contribution based on spatial localization targeted towards problems where the computational cost of classical global methods is prohibitively large.
In this work, we propose a new variant of the trust-region (TR) approach based on the localized reduced basis method (LRBM)  \cite{Ohlberger2015A2865,Ohlberger2018143}. We also refer to \cite{buhr2019localized} for a review of localized model order reduction approaches and to \cite{SmetanaTaddei22} for recent advances in non-linear problems.

In the following, we first review the LRBM for elliptic multiscale problems. We then propose a new TR-LRBM method with certification and adaptive online enrichment for our model 
problem (\ref{P}). Finally, in Section \ref{num_exp}, we show numerical experiments demonstrating the proposed approach's full power. 

\section{Optimality system, primal and dual equation}\vspace*{-.5em}

The primal residual of \eqref{P.state} is key for the optimization as well as for a posteriori error estimation. 
We define for given $u \in V$, $\mu \in \Params$, the primal residual $r_\mu^\pr(u) \in V'$ associated with \eqref{P.state} by
\begin{align}
r_\mu^\pr(u)[v] := l(v;\mu)  - a(u, v;\mu)  &&\text{for all }v \in V.
\label{eq:primal_residual}
\end{align}

Following the approach of \emph{first-optimize-then-discretize}, we base our 
discretization and model order reduction approach on the first-order necessary optimality system, i.e. (cf. \cite{MR4510213} for details and further references)
\begin{subequations}
		\label{eq:optimality_conditions}
		\begin{align}
		r_{\bar \mu}^\pr(\bar u)[v] &= 0 &&\text{for all } v \in V,
		\label{eq:optimality_conditions:u}\\
		\partial_u \J(\bar u,\bar \mu)[v] - a(v,\bar p;\bar \mu)  &= 0 &&\text{for all } v \in V,
		\label{eq:optimality_conditions:p}\\
		(\partial_\mu \J(\bar u,\bar \mu)+\nabla_{\mu} r^\pr_{\bar\mu}(\bar u)[\bar p]) \cdot (\nu-\bar \mu) &\geq 0 &&\text{for all } \nu \in \Params,
		\label{eq:optimality_conditions:mu}
		\end{align}	
	\end{subequations}
where we employ suitable differentiability assumptions. From \eqref{eq:optimality_conditions:p} we deduce the so-called \emph{adjoint equation} and corresponding residual 
\begin{align}
r_\mu^\du(u_\mu, p_\mu)[q] :=  \partial_u \J(u_\mu, \mu)[q] - a(q, p_\mu;\mu) = 0 &&\text{for all } q \in V,
\label{eq:dual}
\end{align}
with solution $p_{\mu} \in V$ for fixed $\mu \in \Params$ and $u_\mu \in V$ solution to \eqref{P.state}.

Since \eqref{P.state} has a unique solution, we use the so-called reduced cost functional $\Jhat: \Params\to \mathbb{R},\,\mu\mapsto \Jhat(\mu) := \J(u_\mu, \mu)$, where we note that the term "reduced" is not associated with model order reduction as introduced below.
For given $\mu\in\Params$, first-order derivatives of~$\Jhat$ can be computed utilizing the adjoint approach, i.e., the gradient $\nabla_\mu\Jhat: \Params \to \R^P$ is given by
	$
	\nabla_{\mu}\Jhat(\mu) = \nabla_{\mu}\J(u_{\mu}, \mu)+\nabla_{\mu}r_\mu^\pr(u_{\mu})[p_{\mu}].
	$

\section{\FOM~and localized reduced basis method}\vspace*{-.5em}
\label{sec:lrbm}

Let us introduce appropriate localized discretization schemes for the 
primal and dual equation based on non-overlapping and non-conforming domain decomposition 
in the context of a Discontinuous Galerkin approach on a coarse grid. We restrict our demonstration to the primal equation as an analogous method for the dual equation is straightforward.

In order to derive a full order model (FOM) for our model reduction approach, we assume that a non-overlapping decomposition of the underlying
domain $\Omega$ is given by a coarse grid $\Grid$ with subdomains $T_j \in \Grid, j=1,\ldots N_H$. Furthermore, each cell $T_j$ is further decomposed
by a local triangulation $\grid(T_j)$ that resolves all fine-scale features of the multiscale problem (\ref{eq:o2}). We then define the global fine-scale partition
$\grid$ as the union of all its local contributions, i.e., $\grid = \bigcup_{j=1}^{N_H} \grid(T_j)$.
For a triangulation $\grid(\omega)$ of some $\omega \subseteq \Omega$, let $H^p\big(\grid(\omega)\big) := \big\{ v \in L^2(\omega) \;\big|\; v|_t \in H^p(t) \;\;\forall t \in \grid(\omega)\big\}$ denote the broken Sobolev space of order $p \in \N$ on $\grid(\omega)$. Then, $H^p(\grid)$ naturally inherits the decomposition
$
H^p(\grid) = \bigoplus_{j=1}^{N_H} H^p(\grid(T_j)).
$
We further define the coarse approximation space  $\V^c := Q^1(\Grid) \subset H^2(\Grid)$ and on each coarse element $T$ the fine-scale correction space $\V^f(T) := S^1(\grid(T)) \subset H^1(T)$, 
where $Q^1(\Grid)$ denotes the standard non-conforming and $S^1(\grid(T))$ the standard conforming piecewise linear (discontinuous) Galerkin finite element space on $\Grid$ and $\grid(T)$, respectively. 
Finally, we define the non-conforming solution space\\[-1em]
$$V(\grid):= \bigoplus_{j=1}^{N_H} V^j \subset H^2(\grid) \quad \text{with }\quad  V^j := \V^c|_{T_j} \oplus \V^f(T_j). $$ 
\begin{definition}[Discontinuous Galerkin multiscale FEM (\FOM)]
We call $u_{h, \mu} \in V(\grid)$ an approximate \FOM~reference solution of (\ref{eq:o2}), if 
\begin{align}\label{eq:broken}
  \bdg\big(u_{h, \mu},v;\mu\big) &= \rhs(v; \mu) &&\text{for all }v \in V(\grid).
\end{align}
\end{definition}
Here, the DG bilinear form $\bdg$ and the right hand side $\rhs$ are given as
\begin{align*}
    \bdg(v,w;\mu) \!:= \!\!&\sum_{t\in\grid}\! \int_t \! A(\mu)  \nabla v\cdot\nabla w
    + \!\!\!\! \sum_{e\in\faces(\grid)} \!\!\!\!\!\! \ \bdg^e(v, w; \mu), \ \
    \rhs(v; \mu) \!:= \!\!\sum_{t\in\grid} \!\! \int_t f(\mu)v,
   \intertext{where $\faces(\cdot)$ denotes the set of all faces of a triangulation and the DG coupling bilinear form $\bdg^e$ for a face $e$ is given by}\nonumber \\[-2em]
   \bdg^e(v, w; \mu) &:= \int_e \big< A(\mu) \nabla v\cdot {\n_e}\big>[w]
    + \big< A(\mu) \nabla w\cdot {\n_e} \big>[v]
    + \frac{\sigma_\face(\mu)}{|e|^{\beta}} [v][w].
\end{align*}
For any triangulation $\grid(\omega)$ of some $\omega \subseteq \Omega$, we assign to each face $e \in \faces\big(\grid(\omega)\big)$ a unique normal $\n_e$ pointing away from the adjacent cell $t^-$, where an inner face is given by $e = t^- \cap t^+$ and a boundary face is given by $e = t^- \cap \partial\omega$, for appropriate cells $t^\pm \in \grid(\omega)$.
The average and jump of a two-valued function $v \in H^2\big(\grid(\omega)\big)$ are given by $\big<v\big> := \tfrac{1}{2}(v|_{t^-} + v|_{t^+})$ and $[v] := v|_{t^-} - v|_{t^+}$ for an inner face and by $\big<v\big> := [v] := v$ for a boundary face, respectively.
The parametric penalty function $\sigma_e(\mu)$ and the parameter $\beta$ must be chosen appropriately to ensure coercivity of $\bdg$ and may involve $A$.
We restrict ourselves to the symmetric interior penalty DG scheme for simplicity. For other variants, we refer to \cite{Ohlberger2015A2865} and the references therein.

Based on the definition of the \FOM~above, the localized reduced basis method constructs in an iterative online 
enrichment procedure appropriate low dimensional
local approximation spaces $\redspace_N^j \subset V^j $ of dimensions $N^j$ that form the global reduced solution 
space via \vspace*{-.5em}
\begin{align}
  \redspace_N = \bigoplus_{j=1}^{N_H} \redspace_N^j, &&\quad N := \dim(\redspace_N) = \sum_{j=1}^{N_H}  N^j.
 \end{align}
Once such a reduced approximation space is constructed, the LRBM approximation is defined as follows.
\begin{definition}[The localized reduced basis method (LRBM)]
We call $u_{N, \mu} \in \fullredspace_N$ a localized reduced basis multiscale approximation of (\ref{eq:broken})
if it holds
\begin{align}\label{eq:reduced}
  \bdg(u_{N, \mu},v_N;\mu) &= \rhs(v_N ; \mu) &&\text{for all } v_N \in \redspace_N.
\end{align}
\end{definition}
Note that $u_{N, \mu}$ solves a globally coupled reduced problem, where all arising quantities can nevertheless be locally computed w.r.t.~the local reduced spaces~$\redspace_N^j$.
For details on the construction of the local approximation spaces, e.g., with a Greedy-based procedure, we refer to \cite{Ohlberger2015A2865}.

The local spaces can be built in an online adaptive procedure \cite{Ohlberger2018143}, where only local patch problems need to be solved without requiring a global solve of the \FOM~method.
More precisely, let a reduced space $V_N$ and the corresponding reduced approximation $u_{N, \mu} \in V_N$ of~\eqref{eq:reduced} for a parameter $\mu  \in \Params$ be given. For enriching the local space $V_N^j$ associated with $T \in \Grid$ at $\mu$, we consider a local oversampling domain 
$O_T := U(T)$, where $U(T)$ denotes a neighborhood of $T$, consisting of an additional layer of coarse neighbouring elements $\bar T \in \Grid$ of $T$.
We further define  $\V(O_{T})$ as the restriction of $V(\grid)$ to $O(T)$
and solve for a local correction $\varphi_T \in V(O_T)$, such that
\begin{align}\label{eq:oversampling_problem}
	\bdg \restrict{O_T}  (u_{N, \mu} + \varphi_T, v;\mu) &= \rhs \restrict{O_T} (v ; \mu) &&\text{for all } v \in V(O_T),
\end{align}
with boundary data from the preceding reduced global solution. The local reduced space $V_N^j$ is then enriched 
with the restriction of $\varphi_T$ to $T$.

In what follows, we do not precompute local approximate spaces but build them iteratively within the optimization routine to solve \eqref{P}.

\section{Relaxed trust-region optimization with reduced models}\vspace*{-.5em}
\label{sec:TR-LRB}

A relaxed adaptive trust-region (TR) method using a localized orthogonal decomposition numerical multiscale method has recently been proposed and extensively studied in \cite{KeilOhlberger22}.
The algorithm proposed in the sequel follows this approach by replacing the FOM and the reduced-order model (ROM) from \cite{KeilOhlberger22} by the 
\FOM~and LRBM introduced in Section \ref{sec:lrbm}.

The overall TR-LRBM algorithm iteratively computes a first-order critical point of problem \eqref{P}.
In the following we denote with $\Jhat_N(\mu) := \J(u_{N, \mu}, \mu)$  the reduced objective functional obtained with LRBM, while $\Jhat_h(\mu) := \J(u_{h, \mu}, \mu)$ denotes the reduced functional obtained with \FOM. We will assume that $\Jhat_N$ admits an a posteriori error estimate of the form 
$$
 |\Jhat_h(\mu) -\Jhat_N(\mu) | \leq \Delta_{\Jhat_N}(\mu).
$$
A derivation of a suitable localizable a posteriori error estimate of this form for the LRBM is beyond the scope of this contribution and is subject to a forthcoming article.

For each outer iteration $k\geq 0$ of the relaxed TR method, we consider a model function $m^{(k)}$ as a cheap local approximation of the cost functional~$\Jhat$ in the relaxed trust-region,
which has radius $\delta^{(k)} + \varepsilon^{(k)}$, where $\delta^{(k)}$ can be characterized by the a posteriori error estimator and  $\varepsilon^{(k)}$ denotes a relaxation parameter from an a priori chosen null sequence, where we assume the existence of $K \in \mathbb{N}$ such that $\varepsilon^{(k)} = 0$ for all $k > K$. 
In our approach, we choose
%\begin{equation*}
$
m^{(k)}(\cdot):= \cJhatn^{(k)}(\mu^{(k)}+\cdot)
$
%\end{equation*}
for $k\geq 0$.
The super-index $(k)$ indicates that we use different localized RB spaces $\fullredspace_{N_k}$ in each iteration.
Thus, we can use $\Delta_{\cJhatn^{(k)}}(\mu)$ for characterizing the trust-region.
In every outer iteration step, we solve for a local solution $\bar s \in \Params$ of the following inner error-aware constrained optimization sub-problem:
\begin{equation}
\label{TRsubprob}
\begin{aligned}
\min_{s\in\Params} \cJhatn^{(k)}(\widetilde{\mu}) \quad \text{ s.t. } \quad  &\frac{\Delta_{\cJhatn^{(k)}}(\widetilde{\mu})}{\cJhatn^{(k)}(\widetilde{\mu})}\leq \delta^{(k)} + \varepsilon^{(k)}, \quad
\widetilde{\mu}:= \mu^{(k)}+s \in\Params \\ & \text{ and } r_{\tilde{\mu}}^\pr(u_{\tilde{\mu}})[v]= 0 \, \text{ for all }  v\in V.
\end{aligned}
\end{equation}
and set $\mu^{(k+1)} := \mu^{(k)} + \bar s$ for the next outer iterate.

In our TR algorithm, we build on the algorithm in \cite{MR4269464}. However, we employ a conforming approach, where the reduced primal and dual spaces coincide and are only constructed from snapshots of the primal equation.
In the sequel, we only summarize the main features of the algorithm and refer to \cite{MR4269464} for more details.
A deeper discussion of suitable localized a posteriori error estimates and more sophisticated construction principles are subject to future work.

As usual in the context of the LRBM, we initialize the local RB spaces $V_N^j$ with the Lagrangian partition of unity w.r.t. $\Grid$, interpolated on the local grids.
Typical for the adaptive TR method, we then initialize the spaces with the starting parameter $\mu^{(0)}$ by using local corrections from \eqref{eq:oversampling_problem}.

For every iteration point $\mu^{(k)}$, we solve \eqref{TRsubprob} with the quasi-Newton projected BFGS algorithm combined with an Armijo-type condition and terminate with a standard reduced first-order critical point (FOC) criterion,
modified with a projection on the parameter space $\Proj_\Params$ to account for constraints on the parameter space.
Additionally, we use a second boundary termination criterion to prevent the subproblem from too many iterations on the boundary of the trust-region. 

After the next iterate $\mu^{(k+1)}$ has been computed, we use the sufficient decrease condition to decide whether to accept the iterate:
\begin{align}
\label{Suff_decrease_condition}
\cJhatn^{(k+1)}(\mu^{(k+1)})\leq \cJhatn^{(k)}(\mu_\text{\rm{AGC}}^{(k)}) + \varepsilon^{(k)}&& \text{ for all } k\in\mathbb{N},
\end{align}
where $\mu_\text{\rm{AGC}}^{(k)}$ denotes the approximate generalized Cauchy point computed with one gradient-descent step of \eqref{TRsubprob}.
Condition \eqref{Suff_decrease_condition} can be cheaply checked with the help of a sufficient and necessary condition.
If $\mu^{(k+1)}$ is accepted, we enrich the local RB spaces $V^j_{\red_k}$ by again solving \eqref{eq:oversampling_problem} for every oversampling domain $O(T_j)$ (only for the primal space) and setting 
$
	V^j_{\red_{k+1}} := \vspan\{V^j_{\red_k}, \varphi_{T_j}|_{T_j} \}.
$
We also refer to~\cite{MR4510213}, where a strategy is proposed to skip an enrichment.

Unlike in the previous works, we do not have immediate access to the FOM gradient since we do not rely on expensive global \FOM~solves for enrichment.
Thus, after enrichment, we evaluate the ROM-type FOC condition
\begin{equation}\label{eq:rom-based-termination}
	\|\mu^{(k+1)}-\Proj_\Params(\mu^{(k+1)}-\nabla_\mu \Jhat_N(\mu^{(k+1)}))\|_2\leq \tau_{\text{\rm{FOC}}}.
\end{equation}
If this condition is not fulfilled at the current outer iterate, we continue the algorithm without FOM computations. If it is fulfilled, we can not reliably terminate the algorithm. Instead, we check the overall convergence of the algorithm with the usual FOM-type FOC condition
\begin{equation}\label{eq:fom-based-termination}
	\|\mu^{(k+1)}-\Proj_\Params(\mu^{(k+1)}-\nabla_\mu \Jhat_h(\mu^{(k+1)}))\|_2\leq \tau_{\text{\rm{FOC}}}.
\end{equation}
If fulfilled, we terminate the algorithm. If not, we enrich the space with the primal- and dual \FOM~approximations (which are available from computing $\nabla_\mu \Jhat_h(\mu^{(k+1)})$), i.e.,
$
V^j_{\red_{k+1}} := \vspan\{V^j_{\red_k}, u_{h,\mu^{(k+1)}}|_{T_j}, p_{h,\mu^{(k+1)}}|_{T_j}\}, \text{ for } j=1, \dots, N_H.$
Then, we continue until \eqref{eq:rom-based-termination} is fulfilled again.

Convergence of the algorithm under suitable assumptions can be proven by using Theorem 3.8 in \cite{MR4510213}; see also \cite{KeilOhlberger22}.
Finally, we summarize the basic TR-LRBM algorithm in Algorithm~\ref{alg:basic_TRRB}.\\[-2em]
\begin{algorithm2e}[!h]
	\KwData{
		initial data $\mu^{(0)}, \delta^{(0)}$, relaxation $(\varepsilon^{(k)})_k$, tolerances $\beta_2$, $\tau_\textnormal{{sub}}$, $\tau_\textnormal{{FOC}}$.
	}
	Initialize LRBM model with PoU and by solving \eqref{eq:oversampling_problem} for $\mu^{(0)}$ and set $k=0$\; 
	\While
	{not ROM-based criterion \eqref{eq:rom-based-termination}}{
		Inner iteration: Compute $\mu^{(k+1)} := \mu^{(k)} + s^{(k)}$ with $s^{(k)}$ solution of~\eqref{TRsubprob}% with termination criteria~\eqref{Termination_crit_sub-problem}
		\label{solve_sub_problem}\;
		\uIf{Sufficient decrease condition~\eqref{Suff_decrease_condition} is fulfilled with relaxation $\varepsilon^{(k)}$\label{suff_dec}} {
			Accept $\mu^{(k+1)}$ and
			enrich the LRBM model by solving \eqref{eq:oversampling_problem} at $\mu^{(k+1)}$\; 
			Possibly enlarge the TR-radius \label{accept}\; 
		}
		\Else {
			Reject $\mu^{(k+1)}$, shrink the TR radius $\delta^{(k)}$ and go to Line \ref{solve_sub_problem}\;		
		}
		Set $k=k+1$\;
	}
	\If{not FOM-based criterion \eqref{eq:fom-based-termination}} {
		Enrich LRBM with global solutions from $\Jhat_h$ in \eqref{eq:fom-based-termination} and go back to Line \ref{solve_sub_problem}\;
	}
	\caption{Basic TR-LRBM algorithm} 
	\label{alg:basic_TRRB}
\end{algorithm2e}\\[-4em]

\section{Numerical experiments}\label{num_exp}\vspace*{-.5em}

\begin{figure}[t]
	\begin{subfigure}[b]{0.35\textwidth}
		\includegraphics[width=\textwidth]{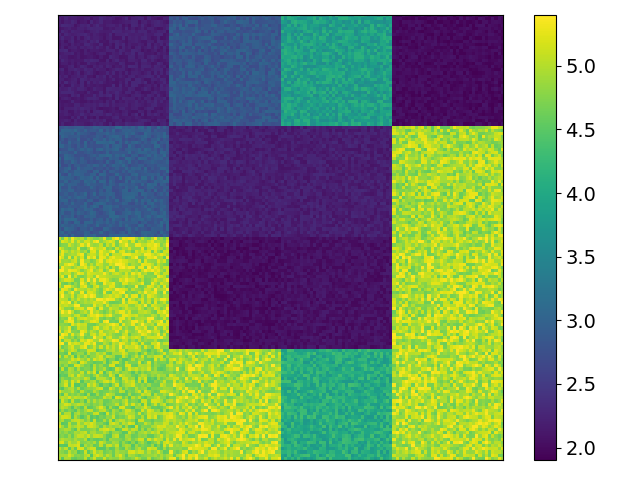}
	\end{subfigure}\qquad
	\begin{subfigure}[b]{0.35\textwidth}
		\includegraphics[width=\textwidth]{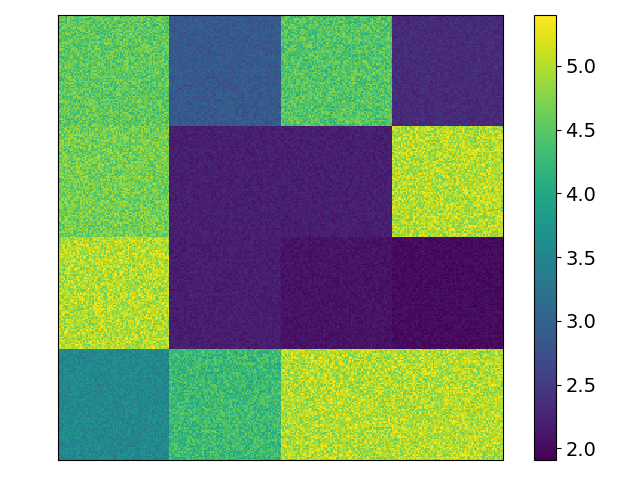}
	\end{subfigure}
	\centering
	\caption{{Coefficient $A^1$ (left) and $A^2$ (right) for the desired state of $\mu^{\text{d}} \in \Params$.}}
	\label{fig:thermal_block}
\end{figure}
We demonstrate the algorithm by using the multiscale benchmark problem used in \cite{KeilOhlberger22}.
The experiments have been conducted with \texttt{PyMOR}\footnote{see \url{https://pymor.org}} \cite{MR3565558} for the model reduction as well as \texttt{dune-gdt}\footnote{see \url{https://github.com/dune-community/dune-gdt}
} and the DUNE framework \cite{bastian2021dune} for the \FOM~discretization.
 We use $\Omega=[0,1]^2$ and define an $L^2$-misfit objective functional with a Tikhonov-regularization term as 
$
	\J(v, \mu) := \frac{\sigma_d}{2} \int_{\Omega}^{} (v - u^{\text{d}})^2 \mathrm{d}x + \frac{1}{2} \sum^{P}_{i=1} \sigma_i (\mu_i-\mu^{\text{d}}_i)^2 + 1,
$
where $\mu^{\text{d}} \in \Params$ denotes the desired parameter with weight~$\sigma_d \in \R$ and $u^{\text{d}} = u_{h}(\mu^{\text{d}})$ the precomputed desired \FOM-FOM solution.
Note that we do not use a coarsening as was used in \cite{KeilOhlberger22}.

For the diffusion coefficient $A$ in \eqref{eq:o2}, we consider a 4x4-thermal block problem with two different thermal block multiscale coefficients $A^1$ and $A^2$, i.e.
\begin{equation*}
%$
A(\mu) := A^{1}(\mu) + A^{2}(\mu) = \sum\limits_{\xi = 1}^{16} \mu_\xi A_\xi^{1} +  \sum_{\xi =17}^{32} \mu_\xi A_{\xi-16}^{2}.
%$
\end{equation*}
We consider $\Params = [1,4]^{24} \times [1,1.2]^8 \subseteq \mathbb{R}^{32}$ and the parameterized multiscale blocks are given by $A_\xi^{1} = A^{1}\big|_{\Omega_{i,j}}$ and $A_\xi^{2} = A^{2}\big|_{\Omega_{i,j}}$, where $\Omega_{i,j}$ denotes the $(i,j)$-th thermal block for $i,j=1,\dots,4$ enumerated by $\xi$.
The multiscale features are uniformly distributed non-periodic values in $[0.9,1.1]$ on $N_1 \times N_1$ (for $A^{1}$) and $N_2 \times N_2$ (for $A^{2}$) quadrilateral grids; see \cref{fig:thermal_block} for a visualization. We set $N_1=150$ and $N_2=300$. Thus, it suffices to have $n_h \times n_h$, $n_h=600$ fine elements for $\grid(\Omega)$ and $n_H \times n_H$, $n_H=10$ coarse elements for $\Grid$.
Both coefficients $A^{1}$ and $A^{2}$ have low-conductivity blocks in the middle of the domain, i.e., for $\Omega_{i,j}$, $i,j=2,3$, which is enforced by a restriction on the parameter space.
Further, we use a constant source $f \equiv 10$.
For all other parameters, we refer to \cite{MR4269464} and its accompanying code \cite{CodeTRLRB}, where the same hyper-parameters for the experiment were used.

We compare two optimization algorithms: An entirely FOM-based BFGS algorithm, where only \FOM~evaluations have been used, and the relaxed TR method described in \cref{sec:TR-LRB} using BFGS for the sub-problems.
For termination, we use $\tau_\text{FOC} = 3 \cdot 10^{-6}$ in order to ensure that the FOM and ROM algorithms stop with the same optimization error in the parameter. 
 \begin{table}[h] \centering \footnotesize
	\begin{tabular}{|l|c|c|c|c|c|c|c|c|c|} %c|
		\hline
		& \multicolumn{3}{c|}{Evaluations} 	& \multicolumn{2}{c|}{Iterations} &   $\mu$-error \\ \hline
		& \FOM~ \eqref{eq:broken}	& LRBM \eqref{eq:reduced}	& Local \eqref{eq:oversampling_problem} 	 & outer & inner &
		\\ \hline  \hline
		Cost factor   & $h$    & $N_{\text{RB}}$  & $O_T$ &  & & \\ \hline
		BFGS with \FOM & 259 & - & - & 85 & -  & 2.89e-3
		\\
		Relaxed TR-LRBM & 2 & 506 & 294 & 2 & 140 & 2.38e-3
		\\ \hline
		%Relaxed TR-RB \cite{MR4269464} & 8 (FEM) & 830 (RB) & - & 3 & 204 & 3.56e-4 \\ \hline
	\end{tabular}
	\caption{Evaluations and accuracy of FOM and ROM.}
	\label{tab}
\end{table}
\begin{figure}[t!]
	\footnotesize
	\centering
	\begin{tikzpicture}
		\definecolor{color0}{rgb}{0.65,0,0.15}
		\definecolor{color1}{rgb}{0.84,0.19,0.15}
		\definecolor{color2}{rgb}{0.96,0.43,0.26}
		\definecolor{color3}{rgb}{0.99,0.68,0.38}
		\definecolor{color4}{rgb}{1,0.88,0.56}
		\definecolor{color5}{rgb}{0.67,0.85,0.91}
		\begin{axis}[
			name=left,
			width=6.5cm,
			height=4cm,
			log basis y={10},
			tick align=outside,
			tick pos=left,
			legend style={nodes={scale=0.7}, fill opacity=0.8, draw opacity=1, text opacity=1,
				xshift=-0.6cm, yshift=-0.1cm, xshift=4.5cm, draw=white!80!black},
			x grid style={white!69.0196078431373!black},
			xlabel={outer iterations},
			xmajorgrids,
			xtick style={color=black},
			y grid style={white!69.0196078431373!black},
			ymajorgrids,
			ymode=log,
			ylabel={\(\displaystyle \| \overline{\mu}-\mu^{(k)} \|^\text{rel}_2\)},
			ytick style={color=black}
			]
			\addplot [semithick, color1, mark=triangle*, mark size=3, mark options={solid, fill opacity=0.5}]
			table {%
				0 7.99350495043297
				1 5.89974734320456
				2 4.78075343104949
				3 4.01547236375528
				4 3.84481197424093
				5 3.78548536642446
				6 3.71970665071376
				7 3.65657943266562
				8 3.61582547732804
				9 3.62798015821625
				10 3.64920505432701
				11 3.66087355519238
				12 3.66096151989762
				13 3.63036213314798
				14 3.51022853521032
				15 3.07239766213454
				16 1.81897274089943
				17 1.36167692213888
				18 1.22000682880227
				19 1.18274249936012
				20 1.05575334568137
				21 0.844509446362641
				22 0.766714779525751
				23 0.742835577447862
				24 0.740871744167076
				25 0.738727508073911
				26 0.734287681134506
				27 0.722043789504619
				28 0.698089485583813
				29 0.659080136514468
				30 0.571739631248157
				31 0.44151526379918
				32 0.389363496376212
				33 0.351576349876455
				34 0.324492117795321
				35 0.310511187976682
				36 0.303359433188223
				37 0.297607372797122
				38 0.286712998458905
				39 0.276419400114934
				40 0.263348266188248
				41 0.239749085651413
				42 0.231119600431893
				43 0.222256775099819
				44 0.210220876537446
				45 0.204419315520867
				46 0.199338673272004
				47 0.192901550179036
				48 0.183909342174311
				49 0.175260224084692
				50 0.1690514357745
				51 0.164871038634135
				52 0.159919490219364
				53 0.14982735240289
				54 0.134031917886848
				55 0.113074587353283
				56 0.100417412296406
				57 0.0916449871441435
				58 0.0862881696266681
				59 0.0804595051714972
				60 0.0740893346782616
				61 0.0681849025569378
				62 0.0641437843286756
				63 0.0584611064319901
				64 0.0462653816344427
				65 0.0348936299001858
				66 0.0196347648264592
				67 0.0121962127730465
				68 0.00894129197717073
			};
			\addlegendentry{BFGS with \FOM}
			\addplot [semithick, black, mark=*, mark size=3, mark options={solid, fill opacity=0.5}]
			table {%
				0 7.99350495043297
				1 0.112462167375487
				2 0.00238527924443378
			};
			\addlegendentry{Relaxed TR-LRBM}
	\end{axis}
	\end{tikzpicture}
	\captionsetup{width=\textwidth}
	\caption{\footnotesize{%
			Relative error decay w.r.t. the optimal parameter of selected algorithms for a single optimization run with the same random initial guess $\mu^{(0)}$ and $\tau_\text{FOC} = 3 \cdot 10^{-6}$.}}
	\label{Fig}
\end{figure}

In \cref{Fig}, we show the number of iterations of the algorithms. Details on the evaluations of global-, reduced- and local problems are given in \cref{tab}.
Both algorithms converged to the same point. The TR algorithm almost entirely avoids FOM evaluations, only required in the FOM-based termination~\eqref{eq:fom-based-termination}, which is only evaluated once (corresponding to $2$ FOM evaluations for primal and dual).

For large-scale problems, where FEM evaluations and even \FOM~evaluations become more and more expensive, we expect our proposed TR method to be compelling. More complex numerical examples and other model problems are subject to further research.
\vspace*{-.5em}

\subsubsection{Acknowledgements} 
The authors acknowledge funding by the BMBF under contracts 05M20PMA and by the Deutsche Forschungsgemeinschaft under contract OH 98/11-1 as well as under Germany’s Excellence Strategy EXC 2044 390685587, Mathematics M\"unster: Dynamics -- Geometry -- Structure.


\end{document}